\documentclass{article}
\usepackage{graphicx,amsmath,amsfonts,amscd,amsthm,amssymb,eucal,pstcol,psfrag,color,url}

\usepackage{tikz}
\usetikzlibrary{matrix}
\tikzset{node style ge/.style={}}
\usepackage{tabstackengine}
\stackMath
\newcommand\w{\makebox[1.3em]{$\color{white}\cdot$}}
\newcommand\ww[1]{\makebox[1.3em]{$#1$}}

\theoremstyle{theorem}

\theoremstyle{definition}

\usepackage{tikz}
\usetikzlibrary{arrows,positioning} 
\tikzset{
    >=stealth',
    punkt/.style={
           rectangle,
           rounded corners,
           draw=black, very thick,
           text width=6.5em,
           minimum height=2em,
           text centered},
    pil/.style={
           ->,
           thick,
           shorten <=2pt,
           shorten >=2pt,}
}

\begin{document}

\title{Farey determinants matrix}
\markright{Farey determinants matrix}
\author{
  R.~Tom\'as\thanks{rogelio.tomas@cern.ch}\\
}

\date{\today}

\maketitle
\begin{abstract}
  A new matrix operation based on inserting columns and rows,
  similarly to the mediant operation between fractions, gives rise to the Farey determinants matrix
  or, equivalently, the matrix of the numerators of the differences of
  Farey fractions.
  This matrix allows to visualize established properties and theorems of 
  Farey fractions allowing for more intuitive
  demonstrations and easier understanding.
  Furthermore it is shown how some Farey determinants matrices contain other lower order
  Farey determinants matrices as block matrices around the main diagonal.
\end{abstract}

\section{The Farey determinants matrix}
The Farey sequence $F_N$ of order $N$  is an ascending sequence of irreducible 
fractions between 0 and 1 whose denominators do not exceed $N$~\cite{hardy}.
Let $h_i/k_i\ <\ h_{i+1}/k_{i+1}$ be two
Farey neighbors then, $h_{i+1}k_i - h_ik_{i+1}=1$.
The next Farey fraction to appear between two Farey neighbors is given by
the mediant as
\begin{equation}\nonumber
 \frac{h_i}{k_i}   < \frac{h_i+h_j}{k_i+k_j} < \frac{h_j}{k_j}\ .
\end{equation}
We define the determinant of any two Farey fractions $h_i/k_i$ and  $h_j/k_j$, both in $F_N$, as
\begin{equation}\nonumber
  d_{ij}(N)=\begin{vmatrix}h_j&h_i\\k_j&k_i\end{vmatrix}= h_jk_i - h_ik_j\ .
\end{equation}
$d_{ij}(N)$ is also the numerator of the difference $h_j/k_j-h_i/k_i$,
so that $d_{ii}=0$ and $d_{i(i\pm1)}=\pm1$.
Note that $(N)$ is omitted from $d(N)$ when it is not necessary.

In the following a matrix operation is introduced
resembling the mediant that allows to iteratively compute $d(N)$
for increasing $N$.
Starting from $F_1=\{0/1,\ 1/1\}$, the skew-symmetric unitary matrix $d(1)$
is given by
\begin{equation}\nonumber
  \begin{matrix}0&-1\\1 & 0\end{matrix}\ ,
\end{equation}
we insert one row and one column in the middle with values given
by the sum of the neighboring horizontal or vertical entries.
This is illustrated as follows,
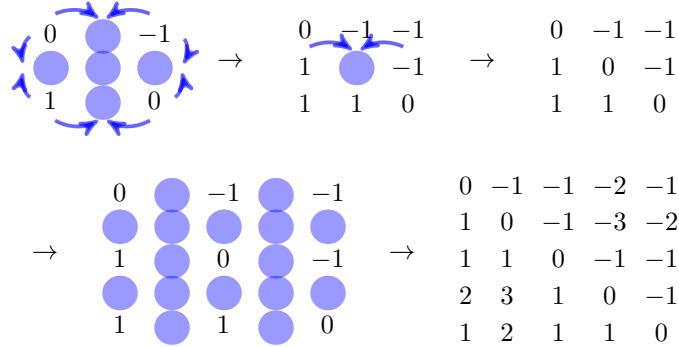
\begin{center}
\begin{tikzpicture}[baseline=(A.center)]
  \tikzset{BarreStyle/.style =   {opacity=.4,line width=4.6 mm,line cap=round,color=#1}}
\matrix (A) [matrix of math nodes, nodes = {node style ge},,column sep=0 mm] 
{  
    0 & \w & -1 \\
 \w &  \w  & \w \\
 1 & \w & 0\\
};
\draw [BarreStyle=blue] (A-2-1.center)  to (A-2-1.center) ;
\draw [BarreStyle=blue] (A-2-2.center)  to (A-2-2.center) ;
\draw [BarreStyle=blue] (A-1-2.center)  to (A-1-2.center) ;
\draw [BarreStyle=blue] (A-3-2.center)  to (A-3-2.center) ;
\draw [BarreStyle=blue] (A-2-3.center)  to (A-2-3.center) ;
\draw (A-1-1.north) edge[pil,bend right=-30,color=blue, line width=0.5mm, opacity=0.6]  (A-1-2.north);
\draw (A-1-3.north) edge[pil,bend right=30,color=blue, line width=0.5mm, opacity=0.6]  (A-1-2.north);
\draw (A-3-1.south) edge[pil,bend right=30,color=blue, line width=0.5mm, opacity=0.6]  (A-3-2.south);
\draw (A-3-3.south) edge[pil,bend right=-30,color=blue, line width=0.5mm, opacity=0.6]  (A-3-2.south);
\draw (A-1-1.west) edge[pil,bend right=0,color=blue, line width=0.5mm, opacity=0.6, out=-40,in=210]  (A-2-1.west);
\draw (A-3-1.west) edge[pil,bend right=0,color=blue, line width=0.5mm, opacity=0.6, out=40,in=-210]  (A-2-1.west);
\draw (A-1-3.east) edge[pil,bend right=0,color=blue, line width=0.5mm, opacity=0.6, out=40,in=-210]  (A-2-3.east);
\draw (A-3-3.east) edge[pil,bend right=0,color=blue, line width=0.5mm, opacity=0.6, out=-40,in=210]  (A-2-3.east);
\end{tikzpicture}
$\ \ \rightarrow \ \ $
\begin{tikzpicture}[baseline=(A.center)]
  \tikzset{BarreStyle/.style =   {opacity=.4,line width=4.6 mm,line cap=round,color=#1}}
\matrix (A) [matrix of math nodes, nodes = {node style ge},,column sep=0 mm] 
{  
 \ww{0} & -1 & -1 \\
    1 &   & -1 \\
    1 & 1 & 0\\
};
\draw [BarreStyle=blue] (A-2-2.center)  to (A-2-2.center) ;
\draw (A-2-1.north) edge[pil,bend right=-30,color=blue, line width=0.5mm, opacity=0.6]  (A-2-2.north);
\draw (A-2-3.north) edge[pil,bend right=30,color=blue, line width=0.5mm, opacity=0.6]  (A-2-2.north);
\end{tikzpicture}
$\ \ \rightarrow \ \ $
\begin{tikzpicture}[baseline=(A.center)]
  \tikzset{BarreStyle/.style =   {opacity=.4,line width=4.6 mm,line cap=round,color=#1}}
\matrix (A) [matrix of math nodes, nodes = {node style ge},,column sep=0 mm] 
{  
    \ww{0} & -1 & -1 \\
   1 &  0  & -1 \\
 1 & 1 & 0\\
};
\end{tikzpicture}
\end{center}

\begin{center}
$ \ \  \rightarrow \ \ $
\begin{tikzpicture}[baseline=(A.center)]
  \tikzset{BarreStyle/.style =   {opacity=.4,line width=4.6 mm,line cap=round,color=#1}}
\matrix (A) [matrix of math nodes, nodes = {node style ge},,column sep=0 mm] 
{  
0 &  \w & -1 & \w   & -1 \\
\w  & \w  &\w   & \w & \w \\
 1&\w   &  0  & \w & -1 \\
  \w&  \w & \w  & \w  &\w  \\
 1& \w  &  1 & \w & 0\\
};
\draw [BarreStyle=blue] (A-1-2.center)  to (A-1-2.center) ;
\draw [BarreStyle=blue] (A-2-2.center)  to (A-2-2.center) ;
\draw [BarreStyle=blue] (A-2-1.center)  to (A-2-1.center) ;
\draw [BarreStyle=blue] (A-4-1.center)  to (A-4-1.center) ;
\draw [BarreStyle=blue] (A-3-2.center)  to (A-3-2.center) ;
\draw [BarreStyle=blue] (A-4-2.center)  to (A-4-2.center) ;
\draw [BarreStyle=blue] (A-5-2.center)  to (A-5-2.center) ;
\draw [BarreStyle=blue] (A-1-4.center)  to (A-1-4.center) ;
\draw [BarreStyle=blue] (A-2-4.center)  to (A-2-4.center) ;
\draw [BarreStyle=blue] (A-3-4.center)  to (A-3-4.center) ;
\draw [BarreStyle=blue] (A-4-4.center)  to (A-4-4.center) ;
\draw [BarreStyle=blue] (A-5-4.center)  to (A-5-4.center) ;
\draw [BarreStyle=blue] (A-2-3.center)  to (A-2-3.center) ;
\draw [BarreStyle=blue] (A-4-3.center)  to (A-4-3.center) ;
\draw [BarreStyle=blue] (A-2-5.center)  to (A-2-5.center) ;
\draw [BarreStyle=blue] (A-4-5.center)  to (A-4-5.center) ;
\end{tikzpicture}
$\ \  \rightarrow \ \ $
\begin{tikzpicture}[baseline=(A.center)]
  \tikzset{BarreStyle/.style =   {opacity=.4,line width=4.6 mm,line cap=round,color=#1}}
\matrix (A) [matrix of math nodes, nodes = {node style ge},,column sep=0 mm] 
{  
      0 & -1 & -1 & -2 & -1 \\
      1 &  0 & -1 & -3 & -2\\
      1 &  1 &  0 & -1 & -1 \\
      2 &  3 &  1 &  0 & -1\\
      1 &  2 &  1 &  1 &  0\\
};
\end{tikzpicture}
\end{center}
where $d(2)$ and $d(3)$ have been generated as the reader can verify. Between $d(N)$ and $d(N+1)$ new rows and columns should be inserted at the same positions as the new
fractions appearing between $F_N$ and $F_{N+1}$.
By construction the top row consists
of the numerators of Farey fractions with opposite sign.
The bottom row consists of numerators of Farey fractions
in reverse order. Similarly happens for the
first and last columns.
The denominators corresponding to the numeratos in the bottom row 
can be obtained by subtracting the top row to the bottom row.
This is illustrated with $d(4)$ together with the Farey fractions corresponding
to the first and last elements in rows and columns,
\begin{equation}\nonumber
  \begin{matrix}
  &  \frac{0}{1} & \frac{1}{4} &\frac{1}{3} &\frac{1}{2} & \frac{2}{3}& \frac{3}{4} &\frac{1}{1} & \\[0.3cm] 
\scriptstyle0/1 &     0 &-1&-1 & -1 & -2 &-3 &-1 & \scriptstyle1/1  \\
\scriptstyle1/4 &     1 & 0 &-1 & -2 & -5 &-8 &-3& \scriptstyle3/4  \\
\scriptstyle1/3 &    1 & 1 & 0 & -1 & -3 &-5 &-2 & \scriptstyle2/3  \\
\scriptstyle1/2 &    1 & 2 & 1 &  0 & -1 & -2 &-1 & \scriptstyle1/2  \\
\scriptstyle2/3 &     2 & 5 & 3 &  1 &  0 &-1 &-1& \scriptstyle1/3  \\
\scriptstyle3/4 &     3 & 8 & 5 &  2 &  1 & 0&-1& \scriptstyle1/4  \\
\scriptstyle1/1 &     1 & 3 & 2 &  1 &  1 & 1& 0& \scriptstyle0/1   \\[0.3cm]
   & \frac{1}{1} & \frac{3}{4} &\frac{2}{3} &\frac{1}{2} & \frac{1}{3}& \frac{1}{4} &\frac{0}{1}
    \end{matrix}
\end{equation}

However it is not yet demonstrated that this process actually generates $d(N)$.
Let $\delta$ be a $n\times n$ matrix built following the above
procedure.
The element $\delta_{ij}$ is a linear combination
of the corresponding top and bottom elements $\delta_{i0}$
and $\delta_{in}$. This linear combination is the same
for all the elements in row $j$, so we can use the elements
$\delta_{0j}$ and $\delta_{nj}$ to reconstruct the linear
combination as follows
\begin{equation}\nonumber
  \delta_{ij} = \delta_{i0}\frac{\delta_{nj}}{\delta_{n0}} +
  \delta_{in}  \frac{\delta_{0j}}{\delta_{0n}}
  = -\delta_{i0}\delta_{nj} + \delta_{in} \delta_{0j}\ ,
\end{equation}
where we have used $\delta_{n0}=-1$ and  $\delta_{0n}=1$, which are true
by construction.
The first row and column correspond
to the numerators of Farey fractions as $\delta_{i0}=-h_i$ and $\delta_{0j}=h_j$, respectively.
The last row and column correspond to the numerators of Farey fractions
in reverse order as $\delta_{in}=h_{n-i+1}=k_i-h_i$ and $\delta_{nj}=-h_{n-j+1}=-(k_j-h_j)$, respectively.
Therefore,
\begin{equation}
  \delta_{ij} = -h_i(k_j-h_j)  +  (k_i-h_i)h_j 
  = -h_i k_j +  k_ih_j =
  \begin{vmatrix}h_j&h_i\\k_j&k_i\end{vmatrix}\ ,
    \nonumber
 \end{equation}
and $\delta_{ij}=d_{ij}$.
By construction $\delta$, or $d$, is
a square, skew-symmetric matrix with rank equal 2
as the new rows are a linear combination
of the neighboring rows for $N>2$.
$d$ is not only skew-symmetric but also
symmetric around the secondary diagonal and therefore
one could keep  only one  quarter of the matrix still being able to
generate higher order $d$'s.
Starting from $d(3)$,
\begin{center}
\begin{tikzpicture}[baseline=(A.center)]
  \tikzset{BarreStyle/.style =   {opacity=0.6,line width=0.6 mm,line cap=round,color=#1}}
\matrix (A) [matrix of math nodes, nodes = {node style ge},,column sep=0 mm] 
{0 & -1 & -1 \\
 1 &  0 & -1\\
 1 &  1 &  0 \\
};
\draw [BarreStyle=blue] (A-1-1.north west)  to (A-2-2.east) ;
\draw [BarreStyle=blue] (A-2-2.east)  to (A-3-1.south west);
\draw [BarreStyle=blue] (A-3-1.south west)  to (A-1-1.north west);
\end{tikzpicture}
\end{center}
we can proceed as before but only taking the quarter of the matrix
highlited with the blue triangle,

\begin{tikzpicture}[baseline=(A.center)]
  \tikzset{BarreStyle/.style =   {opacity=.4,line width=4.6 mm,line cap=round,color=#1}}
\matrix (A) [matrix of math nodes, nodes = {node style ge},,column sep=0 mm] 
{  
 \w & 0 & \w \\
  0 & \ww{1}  & 1 \\
};
\end{tikzpicture}
$\ \ \rightarrow \ \ $
\begin{tikzpicture}[baseline=(A.center)]
  \tikzset{BarreStyle/.style =   {opacity=.4,line width=4.6 mm,line cap=round,color=#1}}
\matrix (A) [matrix of math nodes, nodes = {node style ge},,column sep=0 mm] 
{  
 \w & \w & 0 & \w & \w \\
    & 0 & \w & \w & \w\\
  0 & \w &\ww{1} &\w  & 1 \\
};
\draw [BarreStyle=blue] (A-2-2.center)  to (A-2-2.center) ;
\draw [BarreStyle=blue] (A-3-2.center)  to (A-3-2.center) ;
\draw [BarreStyle=blue] (A-2-3.center)  to (A-2-3.center) ;
\draw [BarreStyle=blue] (A-2-4.center)  to (A-2-4.center) ;
\draw [BarreStyle=blue] (A-3-4.center)  to (A-3-4.center) ;
\draw (A-3-1.south) edge[pil,bend right=30,color=blue, line width=0.5mm, opacity=0.6]  (A-3-2.south);
\draw (A-3-3.south) edge[pil,bend right=-30,color=blue, line width=0.5mm, opacity=0.6]  (A-3-2.south);
\draw (A-3-3.south) edge[pil,bend right=30,color=blue, line width=0.5mm, opacity=0.6]  (A-3-4.south);
\draw (A-3-5.south) edge[pil,bend right=-30,color=blue, line width=0.5mm, opacity=0.6]  (A-3-4.south);
\draw (A-1-3.west) edge[pil,bend right=0,color=blue, line width=0.5mm, opacity=0.6, out=-40,in=210]  (A-2-3.west);
\draw (A-3-3.west) edge[pil,bend right=0,color=blue, line width=0.5mm, opacity=0.6, out=40,in=-210]  (A-2-3.west);
\end{tikzpicture}
$\ \ \rightarrow \ \ $
\begin{tikzpicture}[baseline=(A.center)]
  \tikzset{BarreStyle/.style =   {opacity=.4,line width=4.6 mm,line cap=round,color=#1}}
\matrix (A) [matrix of math nodes, nodes = {node style ge},,column sep=0 mm] 
{  
 \w & \w & 0 & \w & \w \\
    & 0 & 1 & \w & \w\\
  0 & 1 &\ww{1} &2 & 1 \\
};
\draw [BarreStyle=blue] (A-2-4.center)  to (A-2-4.center) ;
\draw (A-3-4.west) edge[pil,bend right=0,color=blue, line width=0.5mm, opacity=0.6, out=40,in=-210]  (A-2-4.west);
\draw (A-2-3.north) edge[pil,bend right=0,color=blue, line width=0.5mm, opacity=0.6, out=40,in=-210]  (A-2-4.north);
\end{tikzpicture}
\begin{equation}
\ \ \rightarrow \ \ \nonumber
\begin{matrix}
   &     & 0     &  \\
   & 0   & 1     & 3  \\
  0& 1   & 1     & 2    & 1
\end{matrix}
\end{equation}

Note that to obtain the 3 we recall
the symmetries of the original matrix and we add the 2 and the 1
below and left of the 3, respectively.

In Section~\ref{index} $d$ is related to the index of
Farey fractions as defined in~\cite{index,generalizedindex}.
Section~\ref{gcd} illustrates an equality among the greatest common divisors
between elements in $d$ as presented in~\cite{partfransums}.  
Section~\ref{maps} shows how some $d(N)$ contain other smaller $d(i)$
thanks to maps in~\cite{arxiv} that preserve the determinant of
two Farey fractions.

\section{$d_{ij}$ and $k$-indexes}\label{index}
The index $\nu(x_i)$ and the  $k$-index $\nu_k(x_i)$ of the $i^{th}$ Farey
fraction in $F_N$
are  introduced in~\cite{index} and~\cite{generalizedindex}, respectively.
We can relate them to the determinant matrix as
\begin{eqnarray}
\nu(x_i)&=&d_{(i-1)(i+1)}\ , \ {\rm for}\ 2 \leq i \leq |F_N|-1,\nonumber\\ 
\nu_k(x_i)&=&d_{(i-1)(i+k-1)}\ ,\ {\rm for}\ 2 \leq i \leq |F_N|-k+1.\nonumber
\end{eqnarray}
Being $\nu_k(x_i)$ a generalized definition of the index: $\nu_2(x_i)=\nu(x_i)$.
The highlighted diagonals in the following $d(3)$ matrix contain part of the $k$-indexes of fractions in $F_3$.
\begin{center}
\begin{tikzpicture}[baseline=(A.center)]
  \tikzset{BarreStyle/.style =   {opacity=.4,line width=4.6 mm,line cap=round,color=#1}}
\matrix (A) [matrix of math nodes, nodes = {node style ge},,column sep=0 mm] 
{0 & -1 & -1 & -2 & -1 \\
 1 &  0 & -1 & -3 & -2\\
 1 &  1 &  0 & -1 & -1 \\
 2 &  3 &  1 &  0 & -1\\
 1 &  2 &  1 &  1 &  0\\
};
\draw [BarreStyle=blue] (A-2-1.center)  to (A-5-4.center) ;
\draw [BarreStyle=green] (A-4-1.center) to (A-5-2.center) ;
\draw [BarreStyle=red] (A-3-1.center) to (A-5-3.center) ;
\draw [BarreStyle=yellow] (A-5-1.center) to (A-5-1.center) ;
\draw [] (A-2-1.west) node[left] {\color{blue}$\nu_1$} ;
\draw [] (A-3-1.west) node[left] {\color{red}$\nu=\nu_2$} ;
\draw [] (A-4-1.west) node[left] {\color{green}$\nu_3$} ;
\draw [] (A-5-1.west) node[left] {\color{yellow}$\nu_4$} ;

\end{tikzpicture}
\end{center}
The sum of the numbers in the red diagonal for $d(N)$ is easily obtained from Theorem~1 in~\cite{index} as
\[
\sum_{i=2}^{|F_N|-1} d_{(i-1)(i+1)} = 3(|F_N|-1)-2N-1\ .
\]
Any two adjacent rows (or columns) of $d$ correspond to
 ordered lists of numerators and denominators of
 Farey neighbours, illustrated as follows.
 \begin{center}
 \begin{tikzpicture}[baseline=(A.center)]
  \tikzset{BarreStyle/.style =   {opacity=.4,line width=4.6 mm,line cap=round,color=#1}}
\matrix (A) [matrix of math nodes, nodes = {node style ge},,column sep=0 mm] 
{0 & -1 & -1 & -2 & -1 \\
 1 &  0 & -1 & -3 & -2\\
 1 &  1 &  0 & -1 & -1 \\
 2 &  3 &  1 &  0 & -1\\
 1 &  2 &  1 &  1 &  0\\
};
\draw [BarreStyle=green] (A-2-1.center) to (A-2-5.center) ;
\draw [BarreStyle=green] (A-3-1.center) to (A-3-5.center) ;
\end{tikzpicture}
$\displaystyle \ \ \rightarrow \ \ \left\{ \frac{1}{1},\frac{0}{1},\frac{-1}{0},\frac{-3}{-1},\frac{-2}{-1}   \right\}$
\end{center}
This is easy to demonstrate by realizing that the first and
last elements of the rows fulfill
\[d_{0i}d_{n(i+1)} - d_{0(i+1)}d_{ni}=1\ ,
\]
as $d_{0i}=h_i$ and $d_{ni}=-(k_i-h_i)$.
Hence $d_{0i}/d_{0(i+1)}$ and $d_{ni}/d_{n(i+1)}$  are Farey neighbors. Therefore, the numbers within the rows
constitute Farey fractions as they are computed as  mediants. Adjacent fractions are also Farey neighbors.
This is equivalent to Lemma~1 in~\cite{generalizedindex} expressed here as
\begin{equation}\nonumber
  \begin{vmatrix}d_{ij}&d_{(i+1)j}\\ d_{i(j+1)}&d_{(i+1)(j+1)}\end{vmatrix}=1\ , {\rm\ for\ }i\leq|F_N|-1\ , \ 
    j\leq|F_N|-1\ .
\end{equation}

\section{Great common divisors among $d_{ij}$}\label{gcd}

According to Theorem 1 in~\cite{partfransums}
the following equality holds between the great common divisors
of  elements $d_{pq}$ with $p$ and $q$ in $\{i > j > k\}$, 
\[ {\rm gcd}\left(d_{kj}, d_{ki}\right) = {\rm gcd}(d_{kj}, d_{ji}) =
{\rm gcd}(d_{ki}, d_{ji})\ .
\]
This property is illustrated using $d(6)$ as follows,
\begin{center}
\begin{tikzpicture}[baseline=(A.center)]
  \tikzset{BarreStyle/.style =   {opacity=0.5,line width=0.5 mm,line cap=round,color=#1}}
\matrix (A) [matrix of math nodes, nodes = {node style ge},,column sep=0 mm] 
{
  0&  -1&  -1&  -1&  -1&  -2&  -1&  -3&  -2&  -3&  -4&  -5&  -1\\
  1&   0&  -1&  -2&  -3&  -7&  -4& -13&  -9& -14& -19& -24&  -5\\
  1&   1&   0&  -1&  -2&  -5&  -3& -10&  -7& -11& -15& -19&  -4\\
  1&   2&   1&   0&  -1&  -3&  -2&  -7&  -5&  -8& -11& -14&  -3\\
  1&   3&   2&   1&   0&  -1&  -1&  -4&  -3&  -5&  -7&  -9&  -2\\
  2&   7&   5&   3&   1&   0&  -1&  -5&  -4&  -7& -10& -13&  -3\\
  1&   4&   3&   2&   1&   1&   0&  -1&  -1&  -2&  -3&  -4&  -1\\
  3&  13&  10&   7&   4&   5&   1&   0&  -1&  -3&  -5&  -7&  -2\\
  2&   9&   7&   5&   3&   4&   1&   1&   0&  -1&  -2&  -3&  -1\\
  3&  14&  11&   8&   5&   7&   2&   3&   1&   0&  -1&  -2&  -1\\
  4&  19&  15&  11&   7&  10&   3&   5&   2&   1&   0&  -1&  -1\\
  5&  24&  19&  14&   9&  13&   4&   7&   3&   2&   1&   0&  -1\\
  1&   5&   4&   3&   2&   3&   1&   2&   1&   1&   1&   1&   0\\
};
\draw [blue,line width=0.5mm,opacity=0.6] (A-10-2) circle (1.5ex);
\draw [red,line width=0.5mm,opacity=0.6] (A-9-2) circle (1.5ex);
\draw [red,line width=0.5mm,opacity=0.6] (A-8-2) circle (1.5ex);
\draw [red,line width=0.5mm,opacity=0.6] (A-7-2) circle (1.5ex);
\draw [red,line width=0.5mm,opacity=0.6] (A-6-2) circle (1.5ex);
\draw [red,line width=0.5mm,opacity=0.6] (A-5-2) circle (1.5ex);
\draw [red,line width=0.5mm,opacity=0.6] (A-4-2) circle (1.5ex);
\draw [red,line width=0.5mm,opacity=0.6] (A-3-2) circle (1.5ex);
\draw [red,line width=0.5mm,opacity=0.6] (A-10-3) circle (1.5ex);
\draw [red,line width=0.5mm,opacity=0.6] (A-10-4) circle (1.5ex);
\draw [red,line width=0.5mm,opacity=0.6] (A-10-5) circle (1.5ex);
\draw [red,line width=0.5mm,opacity=0.6] (A-10-6) circle (1.5ex);
\draw [red,line width=0.5mm,opacity=0.6] (A-10-7) circle (1.5ex);
\draw [red,line width=0.5mm,opacity=0.6] (A-10-8) circle (1.5ex);
\draw [red,line width=0.5mm,opacity=0.6] (A-10-9) circle (1.5ex);
\draw [BarreStyle=blue] (A-9-2.east)  to (A-10-9.north) ;
\draw [BarreStyle=blue] (A-8-2.east)  to (A-10-8.north) ;
\draw [BarreStyle=blue] (A-7-2.east)  to (A-10-7.north) ;
\draw [BarreStyle=blue] (A-6-2.east)  to (A-10-6.north) ;
\draw [BarreStyle=blue] (A-5-2.east)  to (A-10-5.north) ;
\draw [BarreStyle=blue] (A-4-2.east)  to (A-10-4.north) ;
\draw [BarreStyle=blue] (A-3-2.east)  to (A-10-3.north) ;
\end{tikzpicture}
\end{center}
where the number in the blue circle together with another 2 numbers
connected by any blue line define a triplet of numbers for which
the gcd's computed for all possible combinations within the triplet are equal,
e.g.,
\[
{\rm gcd}(14,2)={\rm gcd}(2,8)={\rm gcd}(14,8)=2.
\]

Let the blue circled number be on the antidiagonal $k=n-j+1$, assuming $d$ is
a $n\times n$ matrix. This is illustrated for $d(5)$ as follows,
\begin{center}
\renewcommand\ww[1]{\makebox[0.7em]{$#1$}}
\begin{tikzpicture}[baseline=(A.center)]
  \tikzset{BarreStyle/.style =   {opacity=0.5,line width=0.5 mm,line cap=round,color=#1}}
\matrix (A) [matrix of math nodes, nodes = {node style ge},,column sep=0 mm] 
{
   0  \\
   1 &  0 \\
   1 &  1 &  0\\
   1 &  2 &  1 &  0\\
   2 &  5 &  3 &  1 &  0 \\
   1 &  3 &  2 &  1 &  1 &  0 \\
   3 & 10 &  7 &  4 &  5 &  1 &  0\\ 
   2 &  7 &  5 &  3 &  4 &  1 &  1 &  0\\
   3 & 11 &  8 &  5 &  7 &  2 &  3 &  1 &  0\\
   4 & 15 & 11 &  7 & 10 &  3 &  5 &  2 &  1 &  0  \\
   \ww{1} &  \ww4 &  \ww3 &  \ww2 &  \ww3 &  \ww{1} &  \ww2 &  \ww{1} &  \ww1 &  \ww1 &  \ww0\\
};
\draw [blue,line width=0.5mm,opacity=0.6] (A-10-2) circle (1.5ex);
\draw [red,line width=0.5mm,opacity=0.6] (A-9-2) circle (1.5ex);
\draw [red,line width=0.5mm,opacity=0.6] (A-8-2) circle (1.5ex);
\draw [red,line width=0.5mm,opacity=0.6] (A-7-2) circle (1.5ex);
\draw [red,line width=0.5mm,opacity=0.6] (A-6-2) circle (1.5ex);
\draw [red,line width=0.5mm,opacity=0.6] (A-5-2) circle (1.5ex);
\draw [red,line width=0.5mm,opacity=0.6] (A-4-2) circle (1.5ex);
\draw [red,line width=0.5mm,opacity=0.6] (A-3-2) circle (1.5ex);
\draw [red,line width=0.5mm,opacity=0.6] (A-10-3) circle (1.5ex);
\draw [red,line width=0.5mm,opacity=0.6] (A-10-4) circle (1.5ex);
\draw [red,line width=0.5mm,opacity=0.6] (A-10-5) circle (1.5ex);
\draw [red,line width=0.5mm,opacity=0.6] (A-10-6) circle (1.5ex);
\draw [red,line width=0.5mm,opacity=0.6] (A-10-7) circle (1.5ex);
\draw [red,line width=0.5mm,opacity=0.6] (A-10-8) circle (1.5ex);
\draw [red,line width=0.5mm,opacity=0.6] (A-10-9) circle (1.5ex);
\draw [BarreStyle=blue] (A-9-2.east)  to (A-10-9.north) ;
\draw [BarreStyle=blue] (A-8-2.east)  to (A-10-8.north) ;
\draw [BarreStyle=blue] (A-7-2.east)  to (A-10-7.north) ;
\draw [BarreStyle=blue] (A-6-2.east)  to (A-10-6.north) ;
\draw [BarreStyle=blue] (A-5-2.east)  to (A-10-5.north) ;
\draw [BarreStyle=blue] (A-4-2.east)  to (A-10-4.north) ;
\draw [BarreStyle=blue] (A-3-2.east)  to (A-10-3.north) ;
\end{tikzpicture}
\end{center}

The symmetry of $d$ around
the antidiagonal $d_{(n-j+)i}=d_{j(n-i+1)}$ implies 
\[ {\rm gcd}(d_{(n-j+1)j}, d_{(n-j)i}) = {\rm gcd}(d_{(n-j)j}, d_{(n-j)(n-i+1)}) =
{\rm gcd}(d_{(n-j)i}, d_{(n-j)(n-i)})\ ,
\]
so the previous property can be seen as applied to a single column (or row),
as shown in the following illustration. 
\begin{center}
\renewcommand\ww[1]{\makebox[0.5em]{$#1$}}
\begin{tikzpicture}[baseline=(A.center)]
  \tikzset{BarreStyle/.style =   {opacity=0.5,line width=0.5 mm,line cap=round,color=#1}}
\matrix (A) [matrix of math nodes, nodes = {node style ge},,column sep=0 mm] 
{
   0  \\
   1 &  0 \\
   1 &  1 &  0\\
   1 &  2 &  1 &  0\\
   2 &  5 &  3 &  1 &  0 \\
   1 &  3 &  2 &  1 &  1 &  0 \\
   3 & 10 &  7 &  4 &  5 &  1 &  0\\ 
   2 &  7 &  5 &  3 &  4 &  1 &  1 &  0\\
   3 & 11 &  8 &  5 &  7 &  2 &  3 &  1 &  0\\
   4 & 15 & 11 &  7 & 10 &  3 &  5 &  2 &  1 &  0  \\
   \ww{1} &  \ww4 &  \ww3 &  \ww2 &  \ww3 &  \ww{1} &  \ww2 &  \ww{1} &  \ww1 &  \ww1 &  \ww0\\
};
\draw [blue,line width=0.5mm,opacity=0.6] (A-10-2) circle (1.5ex);
\draw [red,line width=0.5mm,opacity=0.6] (A-9-2) circle (1.5ex);
\draw [red,line width=0.5mm,opacity=0.6] (A-8-2) circle (1.5ex);
\draw [red,line width=0.5mm,opacity=0.6] (A-7-2) circle (1.5ex);
\draw [red,line width=0.5mm,opacity=0.6] (A-6-2) circle (1.5ex);
\draw [red,line width=0.5mm,opacity=0.6] (A-5-2) circle (1.5ex);
\draw [red,line width=0.5mm,opacity=0.6] (A-4-2) circle (1.5ex);
\draw [red,line width=0.5mm,opacity=0.6] (A-3-2) circle (1.5ex);
\draw [blue,line width=0.5mm, opacity=0.6] (A-9-2.east)  to [bend right=80] (A-3-2.east) ;
\draw [blue,line width=0.5mm, opacity=0.6] (A-8-2.east)  to [bend right=80] (A-4-2.east) ;
\draw [blue,line width=0.5mm, opacity=0.6] (A-7-2.east)  to [bend right=50] (A-5-2.east) ;
\path
(A-6-2.east) edge [ -,color=blue, line width=0.5mm, opacity=0.6] node {} (A-6-3.west)  node {} (A-6-2.north);
\end{tikzpicture}
\end{center}

Note that the 3 is paired to himself, and the property above reads:
\[
{\rm gcd}(3,3)={\rm gcd}(15,3)=3\ .
\]
In general, assuming $k=n-j+1$ and $i=(n+1)/2$ implies $(n-i+1) = (n+1)/2$ and
therefore
\[ {\rm gcd}\left(d_{k\, (n-k+1)}, d_{k\, \frac{n+1}{2}}\right) = d_{k\, \frac{n+1}{2}} \ ,
\]
which implies that column-by-column the element in the antidiagonal is a multiple
of the element in the middle row. This last property is also easily shown from the properties of the Farey sequence,
as $d_{k\, (n-k+1)}=b_k(b_k-2a_k)$, with $a_k/b_k$ being the $k^{\rm th}$ element in $F_N$
and $d_{k\, \frac{n+1}{2}}=b_k-2a_k$, since $a_{\frac{n+1}{2}}/b_{\frac{n+1}{2}}=1/2$.

 \section{Maps preserving the determinants matrix}\label{maps}

 Let $F_{N}^{1/a,\, 1/b}$ be the subsequence of $F_N$   defined
as all the fractions of $F_N$ in $[1/a,\ 1/b]$ with $1\leq b \leq a \leq N$.

In~\cite{arxiv} it is demonstrated that
the map
\begin{equation}
F_i \rightarrow F_{N}^{1/q,\, 1/(q-1)}\ , \ \ \ \ \ \ \frac{h}{k} \mapsto \frac{k}{kq-h }\label{map}\ ,
\end{equation}
is bijective between $F_i$ and $F_{N}^{1/q,\, 1/(q-1)}$ when $N$ is a multiple of $i(i+1)$ and $N/(i+1)< q \leq N/i$. It is straight forward to show that this map preserves
the determinant, meaning that
\begin{equation}\nonumber
  \begin{vmatrix}h& h'\\k &k'\end{vmatrix}  =
    \begin{vmatrix}k &k'\\ kq-h  & k'q-h'\end{vmatrix}\ ,
\end{equation}
for $h/k$ and $h'/k'$ belonging to $F_i$.
Therefore $d(N)$
contains $d(i)$ as a matrix block. In other words,
$d(i)$ is contained
$p$ times in $d(i(i+1)p)$.
As an illustration a portion of $d(30)$ is shown containing $d(5)$,
\begin{center}
\begin{tikzpicture}[baseline=(A.center)]
  \tikzset{BarreStyle/.style =   {opacity=0.5,line width=0.5 mm,line cap=round,color=#1}}
\matrix (A) [matrix of math nodes, nodes = {node style ge},,column sep=0 mm] 
{
F_{5} & F_{30} &d  & & & & &  \\          
                 &\scriptstyle4/25 & 0 \\
\scriptstyle0/1  &\scriptstyle1/6  &1 & 0 \\
\scriptstyle1/5 &\scriptstyle5/29 & 9 & 1& 0 \\
\scriptstyle1/4 &\scriptstyle4/23 & 8 &  1  & 1& 0 \\
\scriptstyle1/3 &\scriptstyle3/17 & 7  & 1  & 2 &  1& 0 \\
\scriptstyle2/5 &\scriptstyle5/28 &13 &  2  & 5 &  3 &  1& 0 \\
\scriptstyle1/2 &\scriptstyle2/11 & 6 &  1  & 3 &  2 &  1 &  1& 0 \\
\scriptstyle3/5 &\scriptstyle5/27 &17  & 3 & 10 &  7 &  4 &  5  & 1& 0 \\
\scriptstyle2/3 &\scriptstyle3/16 &11  & 2 &  7  & 5  & 3 &  4 &  1  & 1& 0 \\
\scriptstyle3/4 &\scriptstyle4/21 &16 &  3 & 11  & 8 &  5 &  7 &  2  & 3  & 1& 0 \\
\scriptstyle4/5 &\scriptstyle5/26 &21 &  4 & 15 & 11 &  7 & 10 &  3  & 5  & 2 &  1& 0 \\
\scriptstyle1/1  &\scriptstyle1/5  & 5 &  1  & 4 &  3 &  2 &  3  & 1 &  2 &  1 &  1 &  1& 0 \\
                 &\scriptstyle6/29 &34 &  7 & 29 & 22 & 15 & 23  & 8 & 17 &  9 & 10 & 11 &  1& 0 \\
};
\draw [BarreStyle=blue] (A-2-3.north east)  to (A-13-4.south west) ;
\draw [BarreStyle=blue] (A-13-4.south west) to (A-14-15.north east);
\draw [BarreStyle=blue] (A-2-3.north east)  to (A-14-15.north east);
\end{tikzpicture} 
\end{center}
the left columns show the corresponding Farey fraction in $F_{5}$ and $F_{30}$ according
to the map in Eq.~(\ref{map}) with $q=6$. The blue triangle highlights $d(5)$.


\begin{thebibliography}{20}
\bibitem{hardy} G.H.~Hardy and E.M.~Wright, \emph{An Introduction to the Theory of
 Numbers}, Fifth Edition, Oxford Science Publications, 1996.
 %
\bibitem{index} R. R. Hall \& P. Shiu, \emph{The Index of a Farey Sequence}, Michigan Math. J. {\bf51} (2003).
 %
\bibitem{generalizedindex} A.K. Haynes, \emph{Numerators of differences of nonconsecutive Farey fractions}, International Journal of Number Theory, {\bf6}, No. 3, 05.2010, p. 655-666 (2010).
\bibitem{partfransums} R. Tom\'as, \emph{Partial Franel sums}, 	arXiv:1802.07792 [math.NT] (2018).
  %
\bibitem{arxiv} R. Tom\'as, \emph{Asymptotic behavior of a series of Euler's totient function $\varphi(k)$ times the index of $1/k$ in a Farey sequence},  arXiv:1406.6991v2 [math.NT] (2014). 
%
  \end{thebibliography}
\end{document}